\documentclass[12pt,leqno]{amsart}

\usepackage{amssymb}
 \usepackage{amsthm}

\usepackage{amsmath,amsfonts,amsbsy,latexsym,multicol,fancybox}
\usepackage{graphicx}
\usepackage{epstopdf}
\usepackage{color}
\DeclareGraphicsExtensions{.pdf, .jpg, .pnp, .bmp, .fig}
\usepackage{hyperref}

\usepackage{caption}\usepackage{subcaption}

\newtheorem{theorem}{Theorem}
\newtheorem{assumption}{Assumption}
\newtheorem{definition}{Definition}
\newtheorem{lemma}{Lemma}

\newcommand{\bass}{\begin{assumption}}\newcommand{\eass}{\end{assumption}}
\newcommand{\bde}{\begin{definition}} \newcommand{\ede}{\end{definition}}
\newcommand{\ble}{\begin{lemma}} \newcommand{\ele}{\end{lemma}}
\newcommand{\bth}{\begin{theorem}} \newcommand{\ethe}{\end{theorem}}

\newcommand{\bpf}{\begin{proof}}\newcommand{\epf}{\end{proof}}
\newcommand{\barr}{\begin{array}}\newcommand{\earr}{\end{array}}
\newcommand{\beao}{\begin{eqnarray*}}\newcommand{\eeao}{\end{eqnarray*}\noindent}
\newcommand{\beam}{\begin{eqnarray}}\newcommand{\eeam}{\end{eqnarray}\noindent}
\newcommand{\beqq}{\begin{equation}}\newcommand{\eeqq}{\end{equation}\noindent}

\newcommand{\ov}{\overline} \newcommand{\un}{\underbrace}

\newcommand{\wh}{\widehat}

\newcommand{\ga}{\gamma} 

\newcommand{\D}{\Delta}
  \newcommand{\ep}{\epsilon}

\newcommand{\w}{\omega} \newcommand{\W}{\Omega}

\newcommand{\bfE}{{\mathbb E}}\newcommand{\bbE}{{\mathcal E}} 
\newcommand{\bbf}{{\mathcal F}}

\newcommand{\bbi}{{\mathbb I}}

\newcommand{\bbl}{{\mathcal L}}

 \newcommand{\bbN}{{\mathbb N}}

\newcommand{\bfP}{{\mathbb P}}
                              
 \newcommand{\bbR}{{\mathbb R}}

\graphicspath{{./Images/}}

\begin{document}




\title[Convergence rates of the Semi-Discrete method for SDEs]{Convergence rates of the Semi-Discrete method for stochastic differential equations}
 \author{I. S. Stamatiou}
  \address{University of West Attica, Department of Biomedical Sciences}
 \email{joniou@gmail.com, istamatiou@uniwa.gr}
 \author{N. Halidias}
 \address{University of the Aegean, Department of Statistics and Actuarial-Financial Mathematics}
 \email{nick@aegean.gr}

\begin{abstract}
We study the convergence rates of the semi-discrete (SD) method originally proposed in \emph{Halidias (2012), Semi-discrete approximations for stochastic differential equations and applications, International Journal of Computer Mathematics, 89(6)}.  The SD numerical method was originally designed mainly to reproduce qualitative properties of nonlinear stochastic differential equations (SDEs). The strong convergence property of the SD method has been proved, but except for certain classes of SDEs, the order of the method was not studied. We study the order of $\bbl^2$-convergence and show that it can be arbitrarily close to $1/2.$ The theoretical findings are supported by numerical experiments.
\end{abstract}

\maketitle


\section{Introduction}\label{OSD:sec:intro}
\setcounter{equation}{0}

We are interested in the following class of scalar stochastic differential equations (SDEs),
\beqq  \label{OSD-eq:scalarSDEs}
dx_t =a(t,x_t)dt + b(t,x_t)dW_t, \qquad t\in[0,T],
\eeqq
where $a, b: [0,T]\times\bbR\rightarrow\bbR$ are measurable functions such that (\ref{OSD-eq:scalarSDEs}) has a unique solution and $x_0$ is independent of all $\{W_{t}\}_{t\geq0}.$ SDE (\ref{OSD-eq:scalarSDEs}) has non-autonomous coefficients, i.e. $a(t,x), b(t,x)$ depend explicitly on $t.$
SDEs of the type (\ref{OSD-eq:scalarSDEs}), apart from certain cases, c.f \cite{kloeden_platen:1995}, do not have explicit solutions. Therefore the need for numerical approximations for simulations of the paths $x_t(\w)$ is apparent. We are interested in strong approximations (mean-square) of (\ref{OSD-eq:scalarSDEs}), in the case of nonlinear drift and diffusion coefficients. In the same time we want to reproduce some qualitative properties of the solution process such as domain preservation.

In this direction, we study the semi-discrete (SD) method originally proposed in \cite{halidias:2012} and further investigated  in  \cite{halidias_stamatiou:2016}, \cite{halidias:2014}, \cite{halidias:2015}, \cite{halidias:2015d}, \cite{halidias_stamatiou:2015} and recently in \cite{stamatiou:2018} and \cite{STAMATIOU:2019}. The main idea behind the semi-discrete method is freezing on each subinterval appropriate parts of the drift and diffusion coefficients of the solution at the beginning of the subinterval so as to obtain explicitly solved SDEs. Of course the way of freezing (discretization) is not unique.

The SD method is a fixed-time step explicit numerical method which strongly converges to the exact solution and also preserves the domain of the solution; if for instance the solution process $x_t$ is nonnegative then the approximation process $y_t$ is also nonnegative.

Our main goal is to establish the $\bbl^2$-convergence of the SD method and show that it can be arbitrarily close to $1/2.$

Explicit fixed-step Euler methods fail to strongly converge to solutions of (\ref{OSD-eq:scalarSDEs}) when the drift or diffusion coefficient grows superlinearly \cite[Theorem 1]{hutzenthaler_et_al.:2011}. Tamed Euler methods were proposed to overcome the aforementioned problem, cf. \cite[(4)]{hutzenthaler_jentzen:2015}, \cite[(3.1)]{tretyakov_zhang:2013}, \cite{sabanis:2016} and references therein; nevertheless in general they fail to preserve positivity. We also mention the method presented in \cite{neuenkirch_szpruch:2014} where they use the Lamperti-type transformation to remove the nonlinearity from the diffusion to the drift part of the SDE. Moreover, adaptive time-stepping strategies applied to explicit Euler method are an alternative way to address the problem and there is an ongoing research on that approach, see \cite{fang_giles:2018}, \cite{kelly_lord:2017} and  \cite{kelly_et_al:2018}. Our approach is motivated by the truncated Euler-Maruyama method,  see \cite{mao:2015}, \cite{mao:2016}. At this point, we would like to refer to a different approach in solving stochastic differential equations where the main idea is to reduce, even eliminate in cases, the systematic error that appears in the computation of the mean value of a function of the solution of the SDE, c.f. the recent work \cite{ermakov_pogosian:2019} or \cite{wagner:1988}. 

The outline of the article is the following. In Section \ref{OSD:sec:main} we present the setting and the assumptions, Section \ref{OSD:sec:mainresults} includes among other results our main result, that is Theorem \ref{OSD:theorem:StrongConvergenceOrder}, with their proofs. Section \ref{OSD:sec:numerics} provides a numerical illustration and Section \ref{OSD:sec:conclusion} concluding remarks.

\section{Setting and Assumptions}\label{OSD:sec:main}

Throughout, let $T>0$ and $(\Omega, \bbf, \{\bbf_t\}_{0\leq t\leq T}, \bfP)$ be a complete probability space, meaning that the filtration $ \{\bbf_t\}_{0\leq t\leq T} $ satisfies the usual conditions, i.e. is right continuous and $\bbf_0$ includes all $\bfP$-null sets. Let $W_{t,\w}:[0,T]\times\W\rightarrow\bbR$ be a one-dimensional Wiener process adapted to the filtration $\{\bbf_t\}_{0\leq t\leq  T}.$  Consider SDE  (\ref{OSD-eq:scalarSDEs}), which we rewrite here in its integral form
\beqq\label{OSD-eq:general sde}
x_t=x_0 + \int_{0}^{t}a(s,x_s)ds + \int_{0}^{t}b(s,x_s)dW_s,\quad t\in [0,T],
\eeqq
which admits a unique strong solution. In particular, we assume the existence of a predictable stochastic process $x:[0,T]\times \W\rightarrow \bbR$ such that (\cite[Def. 2.1]{mao:2007}),
$$
\{a(t,x_t)\}\in\bbl^1([0,T];\bbR), \quad \{b(t,x_t)\}\in\bbl^2([0,T];\bbR)
$$
and
$$
\bfP\left[x_t=x_0 + \int_{0}^{t}a(s,x_s)ds + \int_{0}^{t}b(s,x_s)dW_s\right]=1, \quad \hbox{ for every } t\in[0,T].
$$

\bass\label{OSD:assA}
Let $f(s,r,x,y), g(s,r,x,y):[0,T]^2\times\bbR^2\rightarrow\bbR$ be such that $f(s,s,x,x)=a(s,x), g(s,s,x,x)=b(s,x),$ where $f,g$ satisfy the following condition $(\phi\equiv f,g)$
$$
|\phi(s_1,r_1,x_1,y_1) - \phi(s_2,r_2,x_2,y_2)|\leq C_R \Big( |s_1-s_2| + |r_1-r_2| + |x_1-x_2| + |y_1-y_2| \Big)
$$
for any $R>0$ such that $|x_1|\vee|x_2|\vee|y_1|\vee|y_2|\leq R,$ where the quantity $C_R$ depends on $R$ and $x\vee y$ denotes the maximum of $x, y.$
\eass

Let us now recall the SD scheme. Consider the equidistant partition $0=t_0<t_1<...<t_N=T$ and $\D=T/N.$
We assume that for every $n\leq N-1,$ the following SDE
\beqq\label{OSD-eq:SD scheme}
y_t=y_{t_n} + \int_{t_n}^{t} f(t_n, s, y_{t_n}, y_s)ds
+ \int_{t_n}^{t} g(t_n, s, y_{t_n}, y_s)dW_s,\quad t\in(t_n, t_{n+1}],
\eeqq
with $y_0=x_0$ a.s., has a unique strong solution.

In order to compare with the exact solution $x_t,$ which is a continuous time process, we consider the following interpolation process of the semi-discrete approximation, in a compact form,
\beqq\label{OSD-eq:SD scheme compact}
y_t=y_0 + \int_{0}^{t}f(\hat{s}, s,y_{\hat{s}},y_s)ds +
\int_{0}^{t}g(\hat{s},s,y_{\hat{s}},y_s) dW_s,
\eeqq
where $\hat{s}=t_{n}$ when $s\in[t_n,t_{n+1}).$ Process (\ref{OSD-eq:SD scheme compact}) has jumps at nodes $t_n.$ The first and third variable in $f, g$ denote the discretized part of the original SDE. We observe from (\ref{OSD-eq:SD scheme compact}) that in order to solve for $y_t$, we have to solve an SDE and not an algebraic equation, thus in this context, we cannot reproduce implicit schemes, but we can reproduce the Euler scheme if we choose $f(s,r,x,y)=a(s,x)$ and $g(s,r,x,y)=b(s,x).$

In the case of superlinear coefficients the numerical scheme (\ref{OSD-eq:SD scheme compact}) converges to the true solution $x_t$ of SDE (\ref{OSD-eq:general sde}) and this is stated in the following, cf. \cite{halidias_stamatiou:2016},
\bth[Strong convergence]\label{OSD-thm:strong_conv}
	Suppose Assumption \ref{OSD:assA} holds and (\ref{OSD-eq:SD scheme}) has  a unique strong solution for every $n\leq N-1,$ where $x_0\in \bbl^p(\Omega,\bbR).$ Let also
	$$
	\bfE(\sup_{0\leq t\leq T}|x_t|^p) \vee \bfE(\sup_{0\leq t\leq T}|y_t|^p)<A,
	$$
	for some $p>2$ and $A>0.$ Then the semi-discrete numerical scheme (\ref{OSD-eq:SD scheme compact}) converges to the true solution of (\ref{OSD-eq:general sde}) in the $\bbl^2$-sense, that is
	\beqq \label{OSD-eq:strong_conv}
	\lim_{\D\rightarrow0}\bfE\sup_{0\leq t\leq T}|y_t-x_t|^2=0.
	\eeqq
\ethe

Relation (\ref{OSD-eq:strong_conv}) does not reveal the order of convergence. We choose a strictly increasing function $\mu:\bbR_+\rightarrow \bbR_+$ such that 
for every $s,r\leq T$
\beqq \label{OSD-eq:mu}
\sup_{|x|\leq u}\left(|f(s,r,x,y)| \vee |g(s,r,x,y)|\right)\leq \mu(u)(1 + |y|), \qquad u\geq1.
\eeqq

The inverse function of $\mu,$  denoted by $\mu^{-1},$ maps $[\mu(1),\infty)$ to $\bbR_+.$ Moreover, we choose a strictly decreasing function $h:(0,1]\rightarrow[\mu(1),\infty)$ and a constant $\hat{h}\geq 1\vee \mu(1)$ such that
\beqq \label{OSD-eq:h}
\lim_{\D\rightarrow0}h(\D)=\infty \quad \hbox{and}\quad \D^{1/6}h(\D)\leq \hat{h} \quad \hbox{for every}\quad  \D\in(0,1].
\eeqq

Now, we are ready to define the truncated versions of $f, g.$ Let $\D\in(0,1]$ and $f_\D, g_\D$ defined by
\beqq \label{OSD-eq:trunc}
\phi_\D(s,r,x,y):=\phi\left(s,r,(|x|\wedge\mu^{-1}(h(\D)))\frac{x}{|x|},y\right),
\eeqq
for $x,y\in\bbR$ where we set $x/|x|=0$ when $x=0.$

It follows that the truncated functions $f_\D, g_\D$ are bounded in the following way for a given step-size $0<\D\leq1,$
\beam \nonumber
|f_\D(s,r,x,y)| \vee |g_\D(s,r,x,y)|&\leq& \mu(\mu^{-1}(h(\D)))(1 + |y|)\\
\label{OSD-eq:bounded_fg}&\leq &h(\D)(1 + |y|),
\eeam
for all $x,y\in\bbR.$

For the equidistant partition of $[0,T]$ with $\D<1$ consider now the following SDE
\beqq\label{OSD-eq:SD scheme_trunc}
y_t^\D=y_{t_n}^\D + \int_{t_n}^{t} f_\D(t_n, s, y_{t_n}^\D, y_s^\D)ds
+ \int_{t_n}^{t} g_\D(t_n, s, y_{t_n}^\D, y_s^\D)dW_s,\quad t\in(t_n, t_{n+1}],
\eeqq
with $y_0=x_0$ a.s. We assume that (\ref{OSD-eq:SD scheme_trunc}) admits a unique strong solution for every $n\leq N-1$ and rewrite it in compact form,
\beqq\label{OSD-eq:SD scheme compact_trunc}
y_t^\D=y_0 + \int_{0}^{t}f_\D(\hat{s}, s,y_{\hat{s}}^\D,y_s^\D)ds +
\int_{0}^{t}g_\D(\hat{s},s,y_{\hat{s}}^\D,y_s^\D) dW_s.
\eeqq

\bass\label{OSD:assB}
Let the truncated versions $f_\D(s,r,x,y), g_\D(s,r,x,y)$ of $f, g$  satisfy the following condition $(\phi_\D\equiv f_\D,g_\D)$
$$
|\phi_\D(s_1,r_1,x_1,y_1) - \phi_\D(s_2,r_2,x_2,y_2)|\leq h(\D) \Big( |s_1-s_2| + |r_1-r_2| + |x_1-x_2| + |y_1-y_2| \Big)
$$
for all $0<\D\leq 1$ and $x_1, x_2, y_1, y_2\in \bbR,$ where $h(\D)$ is as in (\ref{OSD-eq:h}).
\eass

Let us also assume that the coefficients $a(t,x), b(t,x)$ of the original SDE satisfy the Khasminskii-type condition.

\bass\label{OSD:assC}
We assume the existence of constants $p\geq2$ and $C_K>0$ such that $x_0\in \bbl^p(\Omega,\bbR)$  and
$$
xa(t,x) + \frac{p-1}{2}b(t,x)^2\leq C_K(1 + |x|^2)
$$
for all $(t,x)\in[0,T]\times\bbR$.
\eass

A well-known result follows (see e.g. \cite{mao:2007}) when the SDE (\ref{OSD-eq:general sde}) satisfies the local Lipschitz condition plus the Khasminskii-type condition.

\ble
	Under Assumptions \ref{OSD:assA} (for the coefficients $a(t,x), b(t,x)$) and \ref{OSD:assC} the SDE (\ref{OSD-eq:general sde}) has a unique global solution and for all $T>0,$ there exists a constant $A>0$ such that
	$$
	\sup_{0\leq t\leq T}\bfE |x_t|^p<A.
	$$
\ele

\section{Main results}\label{OSD:sec:mainresults}

In this section we provide the proof of our main result Theorem
\ref{OSD:theorem:StrongConvergenceOrder}. We split the proof is
two steps. First, we prove a general estimate of the error of the
SD method for any $\hat{p}>0.$ Then, we establish the
$\bbl^2$-convergence (\ref{OSD-eq:strong_conv_order}). We denote
the indicator function of a set $A$ by $\bbi_{A}.$ The quantity
$C$ may vary from line to line but it remains independent of the
step-size $\D.$

For ease of notation  in the following  we will avoid the
superscript $\D$ of the approximation process and simply  write
$(y_t).$

Let us define the following stopping time for the solution process
$(y_t^{\D}),$ \beqq  \label{OSD-eq:stopping_timesy}
\rho_{\D,R}=\inf\{t\in [0,T]: |y_t^{\D}|>R \, \hbox { or } \,
|y_{\hat{t}}^{\D}|>R\}. \eeqq

\ble[Error bound for the semi-discrete scheme]
	\label{OSD-lem:Error_SDbound} Let Assumptions \ref{OSD:assA} and \ref{OSD:assB} hold. Let $R>1,$ and
	$\rho_{\D,R}$ as in (\ref{OSD-eq:stopping_timesy}). Then the
	following estimate holds
	$$
	\bfE|y_{s\wedge\rho_{\D,R}}-y_{\wh{s\wedge\rho_{\D,R}}}|^{\hat{p}}
	\leq C(\D^{1/2}h(\D)R)^{\hat{p}}, $$ for any ${\hat{p}}>0,$ where
	$C$ does not depend on $\D.$ \ele

\bpf[Proof of Lemma \ref{OSD-lem:Error_SDbound}] We fix a
${\hat{p}}\geq2.$ Let $n_s$ integer such that
$s\in[t_{n_s},t_{n_s+1}).$ It holds that \beao
&&|y_{s\wedge\rho_{\D,R}}-y_{\wh{s\wedge\rho_{\D,R}}}|^{\hat{p}}=\left| \int_{t_{\wh{n_s\wedge\rho_{\D,R}}}}^{s\wedge\rho_{\D,R}}f_\D(\hat{u},u,y_{\hat{u}},y_u)du + \int_{t_{\wh{n_s\wedge\rho_{\D,R}}}}^{s\wedge\rho_{\D,R}}g_\D(\hat{u},u,y_{\hat{u}},y_u) dW_u\right|^{\hat{p}}\\
&\leq&2^{{\hat{p}}-1}\left|\int_{t_{\wh{n_s\wedge\rho_{\D,R}}}}^{s\wedge\rho_{\D,R}}f_\D(\hat{u},u,y_{\hat{u}},y_u)du\right|^{\hat{p}} + 2^{{\hat{p}}-1}\left|\int_{t_{\wh{n_s\wedge\rho_{\D,R}}}}^{s\wedge\rho_{\D,R}}g_\D(\hat{u},u,y_{\hat{u}},y_u) dW_u\right|^{\hat{p}}\\
&\leq&2^{{\hat{p}}-1}|s\wedge\rho_{\D,R}-t_{\wh{n_s\wedge\rho_{\D,R}}}|^{{\hat{p}}-1}\int_{t_{\wh{n_s\wedge\rho_{\D,R}}}}^{s\wedge\rho_{\D,R}}|f_\D(\hat{u},u,y_{\hat{u}},y_u)|^{\hat{p}}du\\
&& + 2^{{\hat{p}}-1}\left|\int_{t_{\wh{n_s\wedge\rho_{\D,R}}}}^{s\wedge\rho_{\D,R}}g_\D(\hat{u},u,y_{\hat{u}},y_u) dW_u\right|^{\hat{p}}\\
&\leq&C\D^{{\hat{p}}-1}(h(\D))^{\hat{p}}\int_{t_{\wh{n_s\wedge\rho_{\D,R}}}}^{s\wedge\rho_{\D,R}}(1+|y_u|^{\hat{p}})du  +2^{{\hat{p}}-1}\left|\int_{t_{\wh{n_s\wedge\rho_{\D,R}}}}^{s\wedge\rho_{\D,R}}g_\D(\hat{u},u,y_{\hat{u}},y_u)dW_u\right|^{\hat{p}}\\
&\leq&C\D^{{\hat{p}}}(h(\D))^{\hat{p}} +
C\D^{{\hat{p}}}(h(\D)R)^{\hat{p}} +
2^{{\hat{p}}-1}\left|\int_{t_{\wh{n_s\wedge\rho_{\D,R}}}}^{s\wedge\rho_{\D,R}}g_\D(\hat{u},u,y_{\hat{u}},y_u)dW_u\right|^{\hat{p}},
\eeao where we have used the H\"older inequality and the bound
(\ref{OSD-eq:bounded_fg}) for the function $f_\D$. Taking
expectations in the above inequality gives \beao
&&\bfE|y_{s\wedge\rho_{\D,R}}-y_{\wh{s\wedge\rho_{\D,R}}}|^{\hat{p}} \leq C\D^{{\hat{p}}}(h(\D)R)^{\hat{p}}  + 2^{{\hat{p}}-1}\bfE\left|\int_{t_{\wh{n_s\wedge\rho_{\D,R}}}}^{ t_{n_s+1}\wedge\rho_{\D,R} }g_\D(\hat{u},u,y_{\hat{u}},y_u) dW_u\right|^{\hat{p}}\\
&\leq& C\D^{{\hat{p}}}(h(\D)R)^{\hat{p}}  +  2^{{\hat{p}}-1}\un{\left(\frac{{\hat{p}}^{{\hat{p}}+1}}{2({\hat{p}}-1)^{{\hat{p}}-1}}\right)^{{\hat{p}}/2}}_{C_{\hat{p}}}\bfE\left|\int_{t_{\wh{n_s\wedge\rho_{\D,R}}}}^{ t_{n_s+1}\wedge\rho_{\D,R} }|g_\D(\hat{u},u,y_{\hat{u}},,y_u) |^2du\right|^{{\hat{p}}/2}\\
&\leq& C\D^{{\hat{p}}}(h(\D)R)^{\hat{p}}  +  2^{{\hat{p}}-1}C_{\hat{p}}\D^{\frac{{\hat{p}}-2}{2}}\bfE\int_{t_{\wh{n_s\wedge\rho_{\D,R}}}}^{ t_{n_s+1}\wedge\rho_{\D,R} }|g_\D(\hat{u},u,y_{\hat{u}},y_u) |^{\hat{p}} du\\
&\leq& C\D^{{\hat{p}}}(h(\D)R)^{\hat{p}} +
C\D^{{\hat{p}}/2-1}(h(\D))^{\hat{p}}\bfE\int_{t_{\wh{n_s\wedge\rho_{\D,R}}}}^{
	t_{n_s+1}\wedge\rho_{\D,R} }(1 + |y_u |^{\hat{p}}) du\leq
C(\D^{1/2}h(\D)R)^{\hat{p}}, \eeao where in the third step we have
used the Burkholder-Davis-Gundy (BDG)  inequality \cite[Th.
1.7.3]{mao:2007}, \cite[Th. 3.3.28]{karatzas_shreve:1988} on the
diffusion term and in the last step the bound
(\ref{OSD-eq:bounded_fg}) for the function $g_\D.$ Now for
$0<{\hat{p}}<2$ we have that
$$
\bfE|y_{s\wedge\rho_{\D,R}}-y_{\wh{s\wedge\rho_{\D,R}}}|^{\hat{p}}\leq\left(\bfE|y_{s\wedge\rho_{\D,R}}-y_{\wh{s\wedge\rho_{\D,R}}}|^2\right)^{{\hat{p}}/2}\leq
C(\D^{1/2}h(\D)R)^{\hat{p}},
$$
where we have used Jensen inequality for the concave function
$\phi(x)=x^{{\hat{p}}/2}.$ \epf

Let us know provide a moment bound for the approximation process
$(y_t^{\D})$.

\ble[Moment bound for the semi-discrete scheme]
	\label{OSD-lem:moment_SDbound} Let  Assumptions \ref{OSD:assB} and
	\ref{OSD:assC} hold. Then for any $R\leq h(\D)$
	\beqq\label{OSD-eq:approxmombound} \sup_{0\leq \D \leq
		1}\sup_{0\leq t \leq T}\bfE|y^{\D}_t|^{p} \leq C, \eeqq for all
	$T>0.$ \ele

\bpf[Proof of Lemma \ref{OSD-lem:moment_SDbound}] We fix a
$\D\in(0,1]$ and a $T>0.$ Application of the It\^o formula and
(\ref{OSD-eq:SD scheme compact_trunc}) yield

\beao
&&\bfE|y_t|^p \leq \bfE |y_0|^p + \bfE\left(\int_{0}^{t}\left(p|y_s|^{p-1}f_\D(\hat{s}, s,y_{\hat{s}},y_s) + \frac{p(p-1)}{2}|y_s|^{p-2}g_\D^2(\hat{s},s,y_{\hat{s}},y_s)\right) ds\right)\\
&\leq& \bfE |y_0|^p + \bfE\left(\int_{0}^{t}p|y_s|^{p-1}\left(f_\D(\hat{s}, s,y_{\hat{s}},y_s)-f_\D(s, s,y_{s},y_s)+a_\D(s,y_s)\right) ds\right)\\
&& + \bfE\left(\int_{0}^{t}\frac{p(p-1)}{2}|y_s|^{p-2}\left(g_\D(\hat{s},s,y_{\hat{s}},y_s)-g_\D(s,s,y_s,y_s) + b_\D(s,y_s)\right)^2ds\right)\\
&\leq& \bfE |y_0|^p + \bfE\left(\int_{0}^{t}\left(p|y_s|^{p-1} +  \frac{p(p-1)}{2}|y_s|^{p-2}\right)h(\D)(|\hat{s}-s| + |y_{\hat{s}}-y_s|)ds\right)\\
&& + \bfE\left(\int_{0}^{t}p|y_s|^{p-2}\left(y_sa_\D(s,y_s) +
\frac{p-1}{2}b_\D^2(s,y_s)\right)ds\right), \eeao where we have
used Assumption \ref{OSD:assB} and $a_\D, b_\D$ denote the
truncated EM approximations, see \cite{mao:2015}, \cite{mao:2016}.
These functions preserve the Khasminskii-type condition, with a
slightly different constant, see \cite[Lemma 2.4]{mao:2015}.
Bearing this property in mind and using repeatedly the Young
inequality
$$
\alpha^{p-j}\beta\leq \frac{p-j}{p}\alpha^{p} +
\frac{j}{p}\beta^{p/j},
$$
for every $\alpha, \beta\geq0$ and $j=1,2$ we have \beao
&&\bfE|y_t|^p \leq C_1  + C_2 \int_0^t \left(\D h(\D) \bfE|y_s|^{p-1} +  h(\D)\bfE |y_{\hat{s}}-y_s||y_s|^{p-1} +\bfE|y_s|^{p-2}(1 + |y_s|^2)\right)ds\\
&\leq& C_1  + C_2 \int_0^t \sup_{0\leq u\leq s}\bfE|y_u|^{p}ds,
\eeao where we have used (\ref{OSD-eq:h}) and Lemma
\ref{OSD-lem:Error_SDbound} with $R\leq h(\D)$. The inequality
above holds for any $t\in[0,T]$ and the right-hand side in
non-decreasing in $t$ suggesting that \beao
\sup_{0\leq u\leq t}\bfE|y_u^\D|^p &\leq& C_1 + C_2 \int_0^t \sup_{0\leq u\leq s}\bfE|y_u^\D|^{p}ds\\
&\leq& C_1e^{C_2T}\leq C, \eeao by the Gronwall inequality. Since
$C$ is independent of $\D$ inequality
(\ref{OSD-eq:approxmombound}) follows. \epf

\bth[Order of strong
	convergence]\label{OSD:theorem:StrongConvergenceOrder} Suppose
	Assumption \ref{OSD:assB} and Assumption \ref{OSD:assC} hold and
	(\ref{OSD-eq:SD scheme_trunc}) has a unique strong solution for
	every $n\leq N-1,$ where $x_0\in \bbl^p(\Omega,\bbR)$ for some
	$p\geq 14+2\ga.$ Let $\ep\in(0,1/3)$ and define for $\ga>0$
	$$
	\mu(u) = \ov{C}u^{1+\ga}, \quad u\geq0 \quad \hbox { and } \quad
	h(\D)=\ov{C} + \sqrt{\ln \D^{-\ep}}, \quad\D\in(0,1].
	$$
	where $\D\leq1$ and $\hat{h}$ are such that (\ref{OSD-eq:h})
	holds. Then the semi-discrete numerical scheme (\ref{OSD-eq:SD
		scheme compact_trunc}) converges to the true solution of
	(\ref{OSD-eq:general sde}) in the $\bbl^2$-sense with order
	arbitrarily close to $1/2,$ that is \beqq
	\label{OSD-eq:strong_conv_order} \bfE\sup_{0\leq t\leq
		T}|y^\D_t-x_t|^2\leq C\D^{1-\ep}. \eeqq 
\ethe

\bpf[Proof of Theorem \ref{OSD:theorem:StrongConvergenceOrder}]
Denote the difference $\bbE_{t}^{\D}:=y^{\D}_{t}-x_{t}$ and define
the following stopping times \beqq\label{OSD-eq:stopping_timesx}
\tau_R = \inf\{t\in [0,T]: |x_t|>R\}, \quad
\theta_{\D,R}:=\tau_R\wedge\rho_{\D,R}, \eeqq for some $R>1$ big
enough. Let the events $\W$ be defined by $\W_R:=\{\w \in\W:
\sup_{0\leq t\leq T}|x_t|\leq R, \sup_{0\leq t\leq
	T}|y_t^{\D}|\leq R\}.$  We have that \beam \nonumber
&&\bfE\sup_{0\leq t\leq T}|\bbE_t|^2 = \bfE\sup_{0\leq t\leq
	T}|\bbE_t|^2\bbi_{\W_R} +
\bfE\sup_{0\leq t\leq T}|\bbE_t|^2\bbi_{(\W_R)^c}\\
\nonumber&\leq& \bfE\sup_{0\leq t\leq T
}|\bbE_{t\wedge\theta_{\D,R}}|^2 +
\left(\bfE\sup_{0\leq t\leq T}|\bbE_t|^p\right)^{2/p}\left(\bfE(\bbi_{(\W_R)^c})^{2p/(p-2)}\right)^{(p-2)/p}\\
\nonumber&\leq& \bfE\sup_{0\leq t\leq T
}|\bbE_{t\wedge\theta_{\D,R}}|^2 +
\left(\bfE\sup_{0\leq t\leq T}|\bbE_t|^p\right)^{2/p}\left(\bfP(\W_R)^c\right)^{(p-2)/p}\\
\nonumber&\leq& \bfE\sup_{0\leq t\leq T
}|\bbE_{t\wedge\theta_{\D,R}}|^2 +
\left(2^{p-1}\bfE\sup_{0\leq t\leq T}(|y_t|^p + |x_t|^p)\right)^{2/p}\left(\bfP(\W_R)^c\right)^{(p-2)/p}\\
\label{OSD-eq:L2_conv}&\leq& \bfE\sup_{0\leq t\leq T
}|\bbE_{t\wedge\theta_{\D,R}}|^2 + 4\cdot
A^{2/p}\left(\bfP(\W_R)^c\right)^{(p-2)/p}, \eeam where $p>2$ is
as is Assumption \ref{OSD:assC}.  We want to estimate each term of
the right hand side of (\ref{OSD-eq:L2_conv}). It holds that \beao
\bfP(\W_R^c) &\leq& \bfP(\sup_{0\leq t\leq T}|y_t|>R) +  \bfP(\sup_{0\leq t\leq T}|x_t|>R)\\
&\leq& (\bfE\sup_{0\leq t\leq T}|y_t|^k)R^{-k} +  (\bfE\sup_{0\leq
	t\leq T}|x_t|^k)R^{-k}, \eeao for any $k\geq1$ where the first
step comes from the subadditivity of the measure $\bfP$ and the
second step from Markov inequality. Thus for $k=p$ we get
$$
\bfP(\W_R^c)\leq 2AR^{-p}.
$$
We estimate the difference
$|\bbE_{t\wedge\theta_{\D,R}}|^2=|y_{t\wedge\theta_{\D,R}}-x_{t\wedge\theta_{\D,R}}|^2.$
It\^o's formula implies that \beao
&&|\bbE_{t\wedge\theta_{\D,R}}|^2=\int_{0}^{t\wedge\theta_{\D,R}}2|\bbE_s|\left(f_\D(\hat{s},s,y_{\hat{s}},y_s)-f(s,s,x_s,x_s)\right)ds\\
&& +\int_{0}^{t\wedge\theta_{\D,R}}\left( g_\D(\hat{s},s,y_{\hat{s}},y_s)-g(s,s,x_s,x_s) \right)^2 ds\\
&&+ \int_{0}^{t\wedge\theta_{\D,R}}2|\bbE_s|\left( g_\D(\hat{s},s,y_{\hat{s}},y_s)-g(s,s,x_s,x_s)\right)dW_s\\
&\leq&\int_{0}^{t\wedge\theta_{\D,R}}|f_\D(\hat{s},s,y_{\hat{s}},y_s)-f(s,s,x_s,x_s)|^2ds + \int_{0}^{t\wedge\theta_{\D,R}} |\bbE_s|^2 ds  + M_t\\
&&+ \int_{0}^{t\wedge\theta_{\D,R}}
|g_\D(\hat{s},s,y_{\hat{s}},y_s)-g(s,s,x_s,x_s)|^2 ds, \eeao where
$M_t:=2\int_{0}^{t\wedge\theta_{\D,R}}|\bbE_s|\left(
g_\D(\hat{s},s,y_{\hat{s}},y_s)-g(s,s,x_s,x_s)\right)dW_s.$ It
holds that \beao
\bfE \sup_{0\leq t\leq T} |M_t| &\leq& 2\sqrt{32}\cdot\bfE \sqrt{\int_{0}^{T\wedge\theta_{\D,R}}|\bbE_s|^2\left( g_\D(\hat{s},s,y_{\hat{s}},y_s)-g(s,s,x_s,x_s)\right)^2ds}\\
&\leq&\bfE \sqrt{\sup_{0\leq s\leq T} |\bbE_{s\wedge\theta_{\D,R}}|^2 \cdot128\int_{0}^{T\wedge\theta_{\D,R}}\left( g_\D(\hat{s},s,y_{\hat{s}},y_s)-g(s,s,x_s,x_s)\right)^2ds}\\
&\leq&\frac{1}{2}\bfE \sup_{0\leq s\leq T}
|\bbE_{s\wedge\theta_{\D,R}}|^2 +
64\bfE\int_{0}^{T\wedge\theta_{\D,R}}\left(
g_\D(\hat{s},s,y_{\hat{s}},y_s)-g(s,s,x_s,x_s)\right)^2ds, \eeao
thus we get that \beam
&&\nonumber\bfE\sup_{0\leq t \leq T}|\bbE_{t\wedge\theta_{\D,R}}|^2 \leq 2\bfE\sup_{0\leq t\leq T}\int_{0}^{t\wedge\theta_{\D,R}}|f_\D(\hat{s},s,y_{\hat{s}},y_s)-f(s,s,x_s,x_s)|^2ds\\
\nonumber&&+ 130 \cdot \bfE\int_{0}^{T\wedge\theta_{\D,R}} |g_\D(\hat{s},s,y_{\hat{s}},y_s)-g(s,s,x_s,x_s)|^2ds\\
\label{OSD-eq:L2_conv_stoptime}&&
+2\int_{0}^{t\wedge\theta_{\D,R}} \bfE\sup_{0\leq l \leq
	s}|\bbE_l|^2 ds.\eeam Note that
$$
|f_\D(\hat{s},s,y_{\hat{s}},y_s)-f(s,s,x_s,x_s)|^2 =
|f_\D(\hat{s},s,y_{\hat{s}},y_s)-f_\D(s,s,x_s,x_s)+
f_\D(s,s,x_s,x_s)-f(s,s,x_s,x_s)|^2.
$$
If $\mu^{-1}(h(\D))\geq R$ then $f_\D(s,s,x_s,x_s)=f(s,s,x_s,x_s)$
and by Assumption \ref{OSD:assB} we get that
$$
\int_{0}^{t\wedge\theta_{\D,R}}|f_\D(\hat{s},s,y_{\hat{s}},y_s)-f(s,s,x_s,x_s)|^2ds\leq
3h^2(\D)\int_{0}^{t\wedge\theta_R}\Big( |y_s - y_{\hat{s}}|^2 +
|\bbE_s|^2  + |\hat{s}-s|^2\Big)ds
$$
Moreover, it holds that
$$
\int_{0}^{t\wedge\theta_{\D,R}}|\hat{s}-s|^2ds \leq
\sum_{k=0}^{[t/\D-1]}\int_{t_k}^{t_{k+1}\wedge\theta_{\D,R}}
|t_k-s|^2ds.
$$
Taking the supremum over all $t\in[0,T]$ and then expectation we
have \beam\nonumber
&&\bfE\sup_{0\leq t\leq T}\int_{0}^{t\wedge\theta_{\D,R}}|f_\D(\hat{s},s,y_{\hat{s}},y_s)-f(s,s,x_s,x_s)|^2ds\leq 3CTh^2(\D)\D h^2(\D)R^2\\
\nonumber&& + 3h^2(\D)\int_{0}^{T}\bfE\sup_{0\leq l\leq s}|\bbE_{l\wedge\theta_{\D,R}}|^2ds  + 3T\D^2h^2(\D)\\
\label{OSD-eq:L2_conv_stoptime_drift}&\leq& C\D h^4(\D)R^2 + 3h^2(\D)\int_{0}^{T}\bfE\sup_{0\leq l\leq s}|\bbE_{l\wedge\theta_{\D,R}}|^2ds,
\eeam where in the first step we have used Lemma
\ref{OSD-lem:Error_SDbound} for $\hat{p}=2.$ An analogue estimate
of type  (\ref{OSD-eq:L2_conv_stoptime_drift}) holds for the
second integral in (\ref{OSD-eq:L2_conv_stoptime}), that is
\beam\nonumber
&&\bfE\sup_{0\leq t\leq T}\int_{0}^{t\wedge\theta_{\D,R}}|g_\D(\hat{s},s,y_{\hat{s}},y_s)-g(s,s,x_s,x_s)|^2ds\\
\label{OSD-eq:L2_conv_stoptime_diff}&\leq& C\D h^4(\D)R^2 +
3h^2(\D)\int_{0}^{T}\bfE\sup_{0\leq l\leq
	s}|\bbE_{l\wedge\theta_{\D,R}}|^2ds. \eeam Plugging the estimates
(\ref{OSD-eq:L2_conv_stoptime_drift}),
(\ref{OSD-eq:L2_conv_stoptime_diff}) into
(\ref{OSD-eq:L2_conv_stoptime}) gives \beao
\bfE\sup_{0\leq t \leq T}|\bbE_{t\wedge\theta_{\D,R}}|^2 &\leq & C\D h^6(\D) + (132\cdot3 h^2(\D) + 2)\int_{0}^{T}\bfE\sup_{0\leq l\leq s}(\bbE_{l\wedge\theta_{\D,R}})^{2}ds\\
&\leq& C\D h^6(\D)e^{396Th^2(\D)+2T}\leq C\D h^6(\D)e^{h^2(\D)},
\eeao where we have applied the Gronwall inequality and used the
fact that $1<R\leq h(\D).$ Relation (\ref{OSD-eq:L2_conv})
becomes, \beqq\label{OSD-eq:L2_conv_error_est} \bfE\sup_{0\leq
	t\leq T}|\bbE_t|^2 \leq C\D h^6(\D)e^{h^2(\D)}  + CR^{2-p}. \eeqq
Recall that $\mu(u) = \ov{C}u^{1+\ga}$ and $h(\D)=\ov{C} +
\sqrt{\ln \D^{-\ep}},$ for $\ep>0$ to be specified later on. We
bound the first term on the right-hand side of
(\ref{OSD-eq:L2_conv_error_est}) in the following way
$$
C\D h^6(\D)e^{h^2(\D)} \leq C\D (\ln \D^{-\ep})^3\D^{-\ep}\leq
C\D^{1-3\ep},
$$
by choosing $\ep<1/3$, where we used the fact that $0\leq z(\ln
z)^3 \leq z^3$ for big enough $z.$ Moreover, by (\ref{OSD-eq:h})
$$\hat{h}>\D^{1/6}h(\D)>\ov{C}\D^{1/6}>\D^{\frac{(1+\ga)(1-\ep)}{p-2}},$$
whenever $1+\ga<p-2,$ which implies
$$
h(\D)\geq \D^{\frac{(1+\ga)(1-\ep)}{p-2}-\frac{1}{6}}.
$$
By the monotone property of $\mu^{-1}$ we have
$$
\mu^{-1}(h(\D))\geq
\ov{C}^{-\frac{1}{1+\ga}}\D^{\frac{(1-\ep)}{p-2}-\frac{1}{6(1+\ga)}}=R,
$$
for $p$ big enough. Estimate (\ref{OSD-eq:L2_conv_error_est})
becomes \beqq\label{OSD-eq:L2_conv_error_est2} \bfE\sup_{0\leq
	t\leq T}|\bbE_t|^2 \leq C\D^{1-2\ep}  +
C\D^{-(1-\ep)+\frac{p-2}{6(1+\ga)}}. \eeqq Since
$p\geq 14+12\ga$ inequality (\ref{OSD-eq:strong_conv_order}) is true.
\epf

\section{Numerical illustration}\label{OSD:sec:numerics}

We will use the numerical example of \cite[Example 4.7]{mao:2016}, that is we take $a(x)=ax(b-x^2)$ and $b(x)=c x,$ with $a,b,c$ positive and with initial condition  $x_0\in\bbR$ in (\ref{OSD-eq:general sde}), i.e.
\beqq  \label{OSD-eq:exampleSDE}
x_t =x_0 + \int_0^t ax_s(b-x_s^2)ds + \int_0^t cx_sdW_s, \qquad t\geq0.
\eeqq
The above equation, known as the scalar stochastic Ginzburgh-Landau equation, c.f. \cite{kloeden_platen:1995}, has a solution that remains positive (actually there is an explicit solution of $x_t$).

\begin{figure}[ht]
	\centering
	\begin{subfigure}{.45\textwidth}
		\includegraphics[width=1\textwidth]{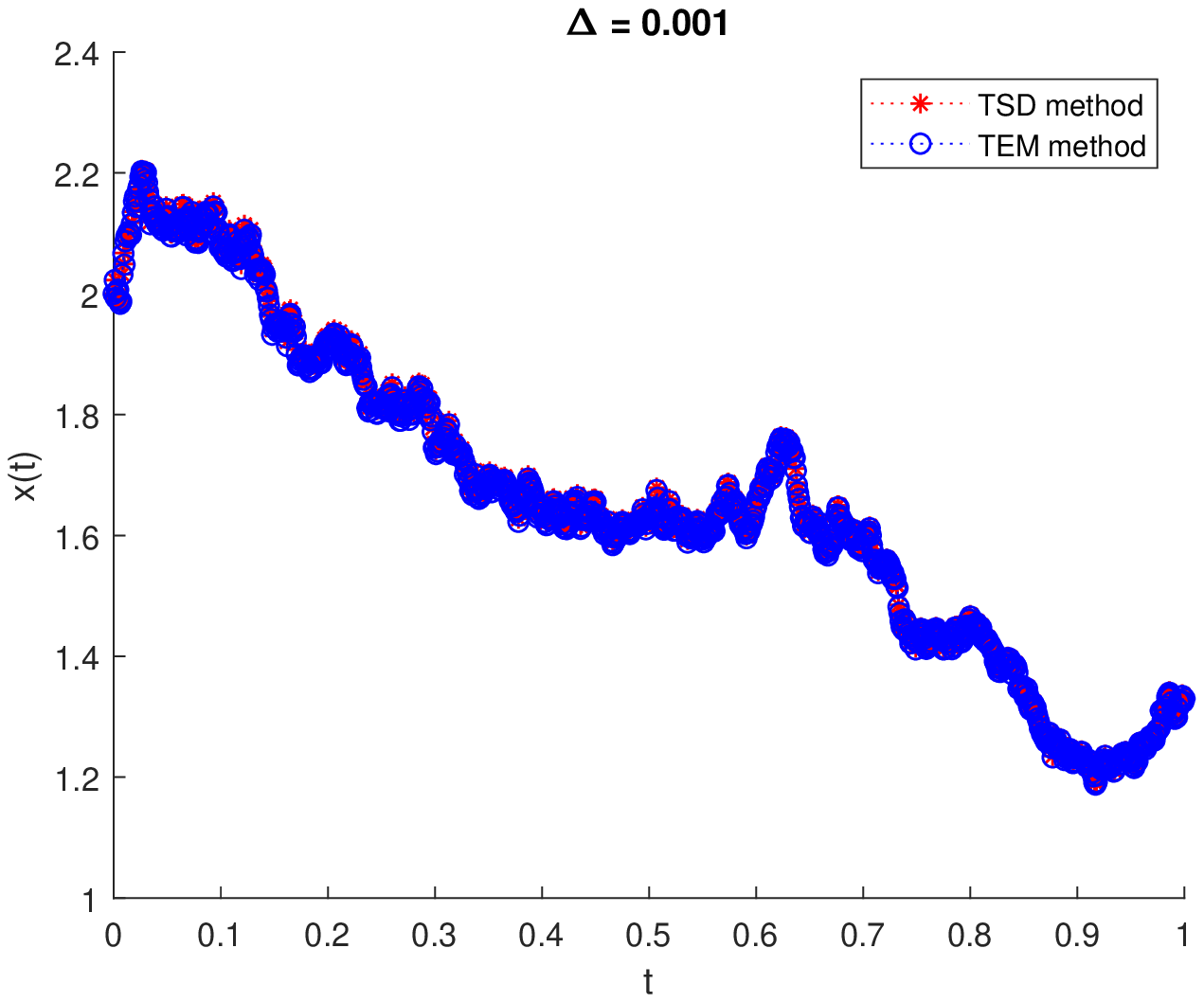}\label{OSD-fig:TSD_TEM1}
		\caption{Trajectory for (\ref{OSD-eq:exampleSDtrunc}) and (\ref{OSD-eq:exampleSDtruncEM}).}
	\end{subfigure}
	\begin{subfigure}{.45\textwidth}
		\includegraphics[width=1\textwidth]{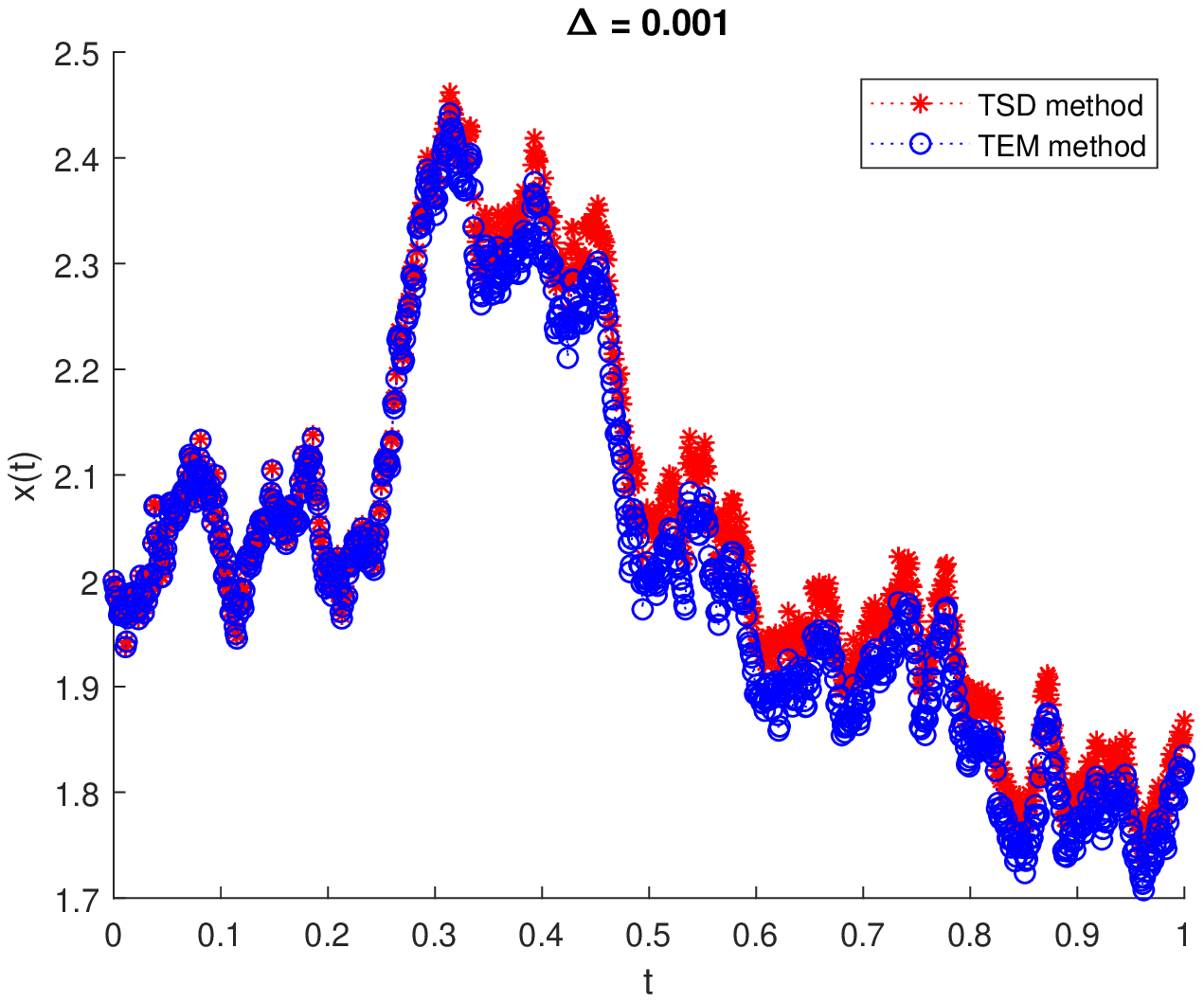}\label{OSD-fig:TSD_TEM2}
		\caption{Trajectory for (\ref{OSD-eq:exampleSDtrunc}) and (\ref{OSD-eq:exampleSDtruncEM}).}
	\end{subfigure}
	\caption{Trajectories of (\ref{OSD-eq:exampleSDtrunc})-(\ref{OSD-eq:exampleSDtruncEM}) for different paths of the Wiener process with $\D=0.001$.}\label{OSD-fig:TSD_TEM}
\end{figure}

Assumption \ref{OSD:assC} holds for any $p>2.$
We choose the auxiliary functions $f,g$ in the following way
$$
f(s,r,x,y) = a(b-x^2)y, \qquad g(s,r,x,y) = cy,
$$
thus (\ref{OSD-eq:SD scheme}) becomes
\beqq\label{OSD-eq:SD schemeExample}
y_t=y_{t_n} + a(b-y_{t_n}^2)\int_{t_n}^{t} y_sds + c\int_{t_n}^{t} y_sdW_s, \quad t\in(t_n, t_{n+1}],
\eeqq
with $y_0=x_0$ a.s., which admits an exponential unique strong solution. In particular,
\beqq  \label{OSD-eq:exampleSD}
y_{n+1} =y_n\exp\left\{\left(a(b-y_{n}^2)-\frac{c^2}{2}\right)\D + c \D W_n\right\}, \quad n\in\bbN,
\eeqq
Note that (\ref{OSD-eq:mu}) holds with $\mu(u) =(a(b+1)\vee c)|u|^3$ since
$$
\sup_{|x|\leq u}\left(|a(b-x^2)y| \vee |cy|\right)\leq (a(b+1)\vee c)|u|^2(1 + |y|), \qquad u\geq1.
$$
Therefore, in the notation of Theorem \ref{OSD:theorem:StrongConvergenceOrder}, $\ga=2$ and $\ov{C}=(a(b+1)\vee c).$ Finally, $h(\D) = \ov{C} + \sqrt{\ln \D^{-\ep_1}}$ for any $\D\in(0,1].$ Clearly $h(1)\geq \mu(1)$ and
$$
\D^{1/6}h(\D) \leq \D^{1/6}\ov{C} + \sqrt{\D^{1/3}\ln \D^{-\ep_1}}\leq \ov{C} + \sqrt{\D^{1/3-\ep_1}}\leq \ov{C} + \D^{1/6-\ep_1/2}\leq  \ov{C} + 1
$$
for any $\D\in(0,1]$ and $0<\ep_1\leq 1/3.$ Therefore we take $\hat{h}=\ov{C} + 1.$
The truncated versions of the semi-discrete method (TSD) read,
\beqq  \label{OSD-eq:exampleSDtrunc}
y_{n+1}^\D =y_n^\D\exp\left\{\left(a(b-\left(y_{n}^\D\wedge \mu^{-1}(h(\D))\right)^2)-\frac{c^2}{2}\right)\D + c \D W_n\right\},
\eeqq
for $n\in\bbN.$
We perform computer simulations for the case $a=0.1, b=1, c=0.2$ and $x_0=2$ as in \cite[Example 4.7]{mao:2016} with $\ep_1 = 1/3$ and compare with the truncated Euler Maruyama method (TEM), which reads
\beqq  \label{OSD-eq:exampleSDtruncEM}
y_{n+1}^{TEM} =y_n + a\left(|y_{n}|\wedge \mu^{-1}(\bar{h}(\D))\frac{y_n}{|y_n|}\right)\D + b\left(|y_{n}|\wedge \mu^{-1}(\bar{h}(\D))\frac{y_n}{|y_n|}\right) \D W_n,
\eeqq
for $n\in\bbN,$ where $\bar{h}(\D)=\D^{-\ep_2/2}$ with $\ep_2=1/2,$ and $\bar{\D}^*\leq (8\ov{C})^{-\frac{2}{\ep_2}}.$
Figure \ref{OSD-fig:TSD_TEM} shows sample simulations paths of $x(t)$ by TSD and TEM respectively with sample size $\D=10^{-3}.$ Note that TSD works for all  $\D< 1$ and TEM works for $\D\leq  0.1526$ as proved in \cite{mao:2016}. (in an updated version of TEM in \cite{hu_li_mao:2018} it is shown that it works for all $\D<1$)

\begin{figure}[ht]
	\centering
	\begin{subfigure}{.45\textwidth}
		\includegraphics[width=1\textwidth]{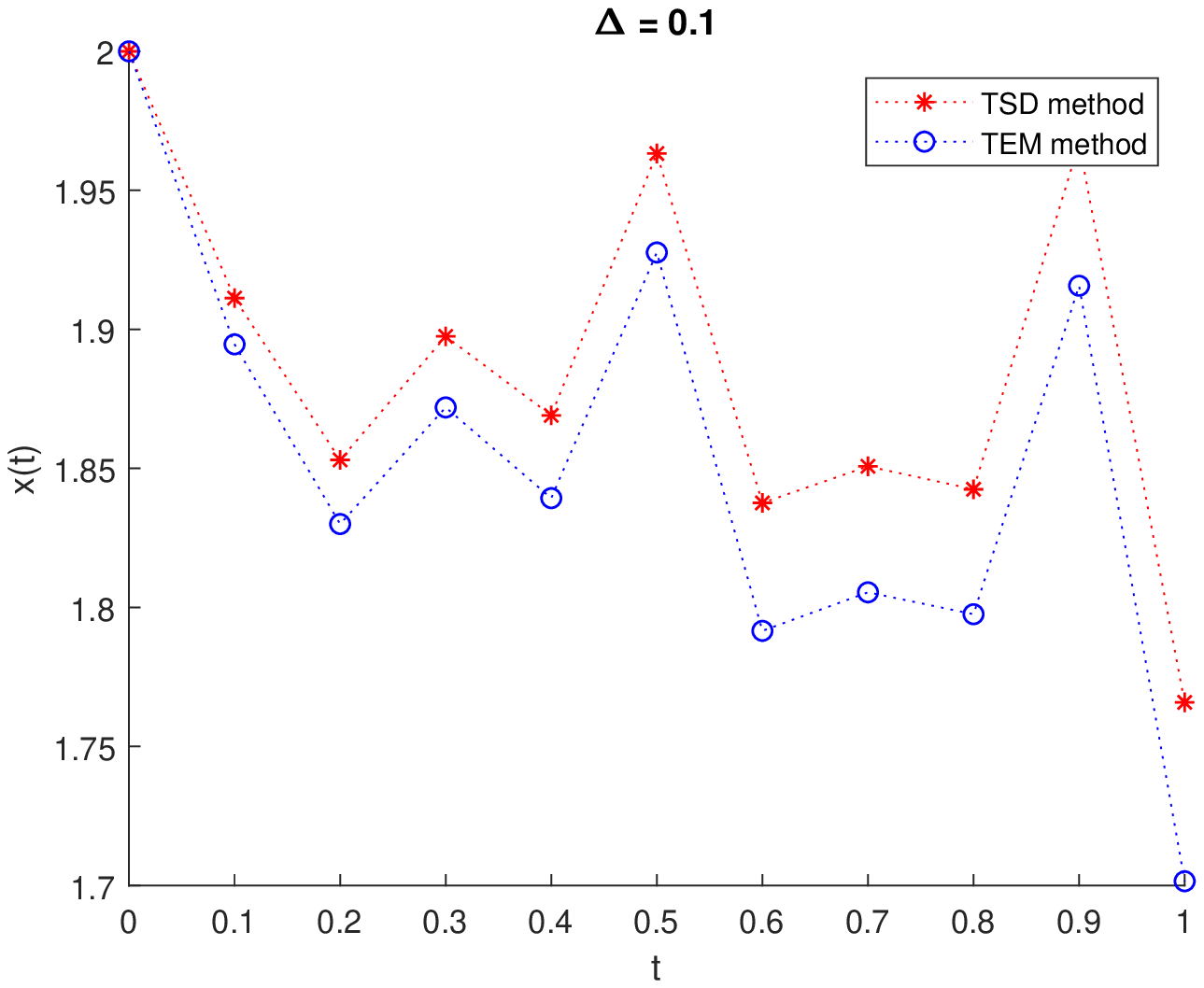}\label{OSD-fig:TSD_TEM3}
		\caption{Trajectory for (\ref{OSD-eq:exampleSDtrunc}) and (\ref{OSD-eq:exampleSDtruncEM}).}
	\end{subfigure}
	\begin{subfigure}{.45\textwidth}
		\includegraphics[width=1\textwidth]{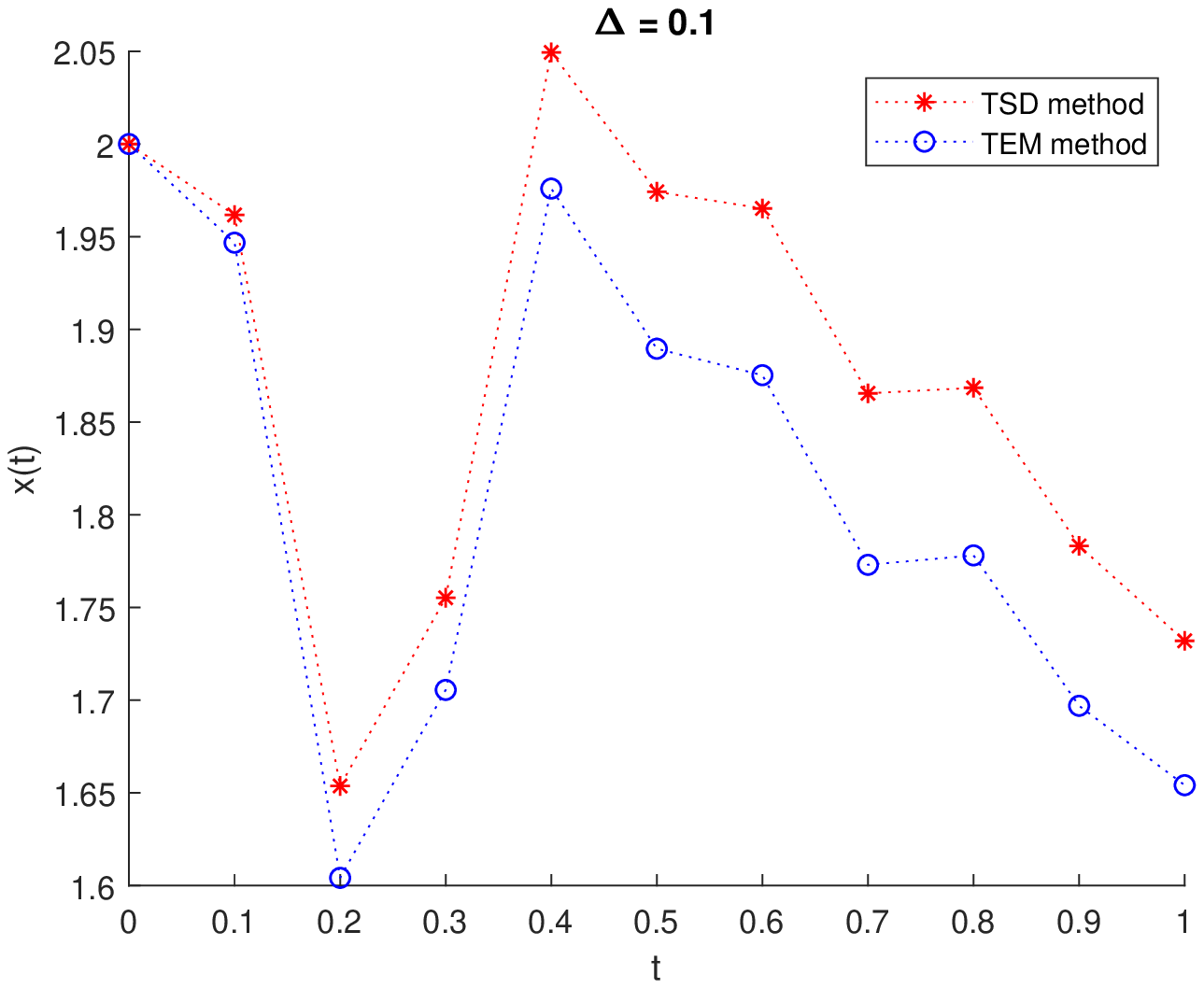}\label{OSD-fig:TSD_TEM4}
		\caption{Trajectory for (\ref{OSD-eq:exampleSDtrunc}) and (\ref{OSD-eq:exampleSDtruncEM}).}
	\end{subfigure}
	\caption{Trajectories of (\ref{OSD-eq:exampleSDtrunc})-(\ref{OSD-eq:exampleSDtruncEM}) for different paths of the Wiener process with $\D=0.1$.}\label{OSD-fig:TSD_TEMd}
\end{figure}

\begin{figure}[ht]
	\centering
	\includegraphics[width=1\textwidth]{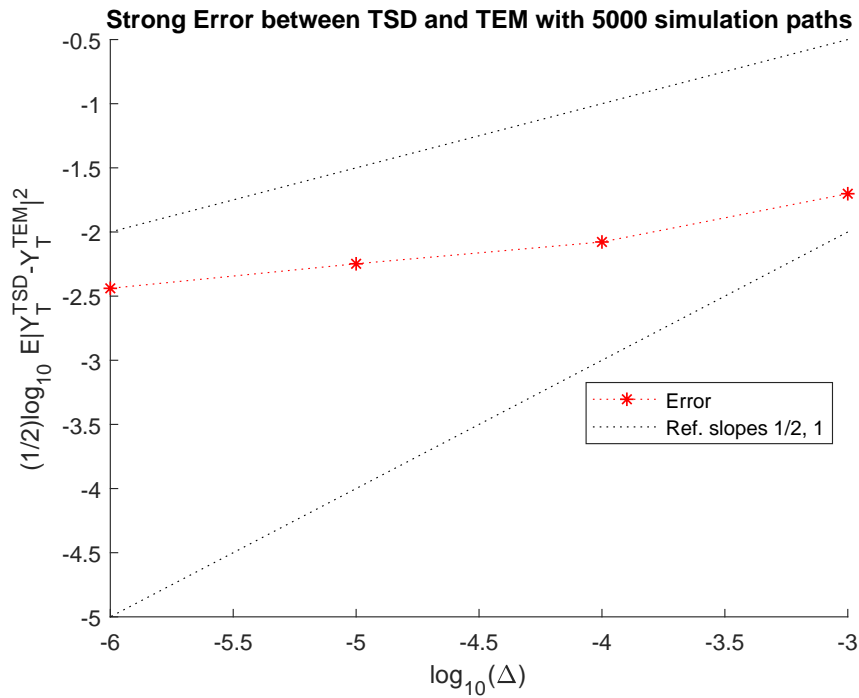}
	\caption{The strong errors between TSD and TEM.}\label{OSD-fig:TSD_TEMerr}
\end{figure}

We also perform $10000$ sample paths of the TSD and TEM respectively for stepsizes $10^{-3}, 10^{-4}, 10^{-5}$ and  $10^{-6}.$ Figure \ref{OSD-fig:TSD_TEMerr} shows the log-log plot of the strong errors between TSD and TEM which is close to $1$ TSD has order $1/2$ in $\bbl^2$-sense thus our TSD has the order $1/2$ in $\bbl^2$-sense too. Nevertheless, the approximation process TEM (\ref{OSD-eq:exampleSDtruncEM}) does not always produce positive values, while TSD (\ref{OSD-eq:exampleSDtrunc}) is positive by construction.

\section{Conclusion and Future Work}\label{OSD:sec:conclusion}
In this paper we study the convergence rates of the semi-discrete (SD) method, originally proposed in \cite{halidias:2012}. Using a truncated version of the  SD method, we show that the order of $\bbl^2$-convergence can be arbitrarily close to $1/2.$ The advantage of our method, over other useful numerical methods (such as the tamed Euler method, the implicit Euler method, the truncated Euler method) applied to nonlinear problems, is that it can reproduce qualitative properties of the solution process. The main qualitative property that has been investigated in all the works so far concerning the SD method is the domain preservation of the solution process. In a future work, we aim to study  other qualitative properties relevant with the stability of the method and answer questions of the following type: \textit{Is the SD method able to preserve the asymptotic stability of the underlying SDE?}

\bibliographystyle{unsrt}\baselineskip12pt
\bibliography{order_of_SD}

\end{document}